\documentclass[A4]{article} 
\usepackage{color,mathrsfs}
\usepackage{amssymb,amsmath,latexsym}
\usepackage{graphicx,amssymb}
\usepackage[utf8]{inputenc}
\usepackage{authblk}
\usepackage[T1]{fontenc}

\newtheorem{thm}{Theorem}
\newtheorem{prop}{Proposition}
\newtheorem{cor}{Corollary}

\newcommand{\N}{\mathbb N}
\newcommand{\R}{\mathbb R}
\newcommand{\Z}{\mathbb Z}

\def\tilde{\widetilde}
\def\epsilon{\varepsilon}

\def\boxit#1{\vbox{\hrule\hbox{\vrule\kern.75truemm
\vbox{\kern.75truemm#1\kern1truemm}\kern1truemm\vrule}\hrule}}

\def\qed{\hskip 1mm\boxit{}\hskip 1mm}

\title{A characterization of the convergence\\ in variation  for the \\ generalized sampling series}

\author{ {\bf Laura Angeloni, \hskip0.3cm Danilo Costarelli, \hskip0.3cm Gianluca Vinti} \\
Department of Mathematics and Computer Science\\
University of Perugia, Via Vanvitelli 1, 06123, Perugia, Italy\\
{\small {\tt laura.angeloni@unipg.it} \hskip0.3cm - \hskip0.3cm {\tt danilo.costarelli@unipg.it}}\\
{\small {\tt gianluca.vinti@unipg.it}}}

\date{}

\begin{document}
\maketitle

\begin{abstract}
In this paper, we study the convergence in variation for the generalized sampling operators based upon averaged-type kernels and we obtain a characterization of absolutely continuous functions. This result is proved exploiting a relation between the first derivative of the above operator acting on $f$ and the sampling Kantorovich series of $f'$. By such approach, also a variation detracting-type property is established. Finally, examples of averaged kernels are provided, such as the central B-splines of order $n$ (duration limited functions) or other families of kernels generated by the Fej\'er and the Bochner-Riesz kernels (bandlimited functions). 
\end{abstract}

\section{Introduction}\label{sec-intro}

The generalized sampling series are an important and well-known family of operators in approximation theory and play a relevant role in Signal Processing. They were introduced around 1980s, when the German mathematician P.L. Butzer established an approximate sampling formula, with the aim to reconstruct not-necessarily bandlimited signals (see e.g., \cite{RIST1,BUFIST}), and are defined as
$$
(S_w f)(t):=\sum_{k\in\Z} f\left({k\over w}\right)\chi(wt-k), \ \ t\in\R,\ w>0.
$$
Here, in place of the sinc-function, as in the classical sampling formula, the operators $S_w$ are based upon the generalized kernel functions $\chi$ which satisfy the classical assumptions of approximate identities (\cite{BUNE}).

For such family of operators, several approximation results were given by means of different kinds of convergence, such as pointwise and uniform convergence, $L^p$ convergence, modular convergence and so on (see \cite{BBSV-06}).   

In this paper we face the problem of the convergence in variation, in the sense of Jordan, for the family of the generalized sampling series and we obtain a complete characterization of the absolutely continuous functions in terms of convergence in variation by means of such family of discrete operators. 

Results about variation of the generalized sampling series were previously studied in some particular cases: for example, the variation detracting property was obtained for some kind of bandlimited kernel functions (\cite{KT-2009,K-2011,KM-2014,ORTA1}), but the topic was never faced in a general setting.

In this paper we consider a general class of generalized sampling series
$$
(\bar S^m_w f)(t):=\sum_{k\in\Z} f\left({k\over w}\right)\bar\chi_m(wt-k), \ \ t\in\R,\ w>0,
$$
based on a family of kernel functions of averaged type, i.e.,
$$
\bar\chi_m(t):={1\over m}  \int_{-{m\over 2}}^{m\over 2} \chi(t+v)\,dv, \ \ t \in \R, \ \ m \in \N,
$$
where $\chi$ satisfies the usual assumptions on kernels and $f \in BV(\R)$ (the space of bounded variation functions on $\R$). 

We first prove (see Section \ref{sec-mres}) that such operators satisfy a variation detracting-type property, namely
$$
V[\bar S^m_w f]\le {1\over m} \Vert \chi\Vert_1 V[f],
$$
for every $w>0$, $m\in\N$. Then, in order to prove the convergence in variation, we first establish a relation between the first derivative of $\bar S^m_w f$ and $K_w f'$, where the operators $K_w$ are the so-called sampling Kantorovich series (see (\ref{sKant}) of Section \ref{sec-not}) based on the kernel $\chi$ (see \cite{BBSV-06,COVI1}). Here the averaged form of the considered kernels $\bar \chi_m$ plays an important role and is not restrictive, as shown below. The sampling Kantorovich operators represent an $L^1$-version of $S_w$, and their approximation properties have been widely studied in last years, both from the theoretical and the applications point of view; see, e.g., \cite{COVI2,ORTA1,COMIVI1}.
  
  By means of the relation between the generalized sampling series and the sampling Kantorovich series, we are able to prove (see Section \ref{sec-mres}) that, for every fixed $m\in\N$, 
$$
V[\bar S_w^m f-f]\rightarrow 0,\ \ w\rightarrow +\infty,
$$
if and only if $f$ belongs to the space $AC(\R)$ of the absolutely continuous functions on $\R$. Note that, for the converse result we use the closedness of $AC(\R)$ in $BV(\R)$ with respect to the variation functional (\cite{BBST-03}), together with the absolute continuity of the generalized sampling series (Proposition \ref{var-dim}).
 
  One of the main advantages that can be reached by the approach proposed in this paper, is the possibility to obtain approximation results also for not-necessarily bandlimited kernels, therefore enlarging the class of kernels (e.g., duration limited kernels). Moreover, we obtain a complete characterization of the space $AC(\R)$ by means of the convergence in variation for $\bar S_w^m$, while previous results only established convergence in variation in some proper subspaces of $AC(\R)$ (i.e., in Bernstein spaces, \cite{BV}).
  
  It is important to mark out that the use of averaged kernels for the operators $S_w$ is not restrictive: indeed there are many examples of kernels widely used in approximation theory that are of averaged type. Among them, for instance, the central B-splines of order $n \in \N$, which are typical examples of duration limited kernels (and therefore not bandlimited), are averaged central B-splines of order $n-1$ (see Section \ref{sec-ex}).  Finally, in Section \ref{sec-ex} we show that many other examples of averaged kernels can be generated by using classical families of kernels, such as the Fej\'er and the Bochner-Riesz kernels, and many others.


\section{Notations and preliminaries}\label{sec-not}

We will work in the frame of the space of functions of bounded variation on $\R$ (see \cite{KR-2011,AD-2014,AV-2014,LS-2015}), namely
$$
BV(\R):=\{f:\R\longrightarrow \R:\ V[f]<+\infty\}.
$$ 
Here $V[f]:=\sup_{[a,b]\subset \R}V_{[a,b]}[f]$ is the Jordan variation of $f$ over $\R$; 
$  V_{[a,b]}[f]=\sup \sum_{i=1}^n |f(x_i)-f(x_{i-1})|$, where the supremum is taken over all the possible partitions $a=x_0<x_1<\ldots<x_n=b$ of the interval $[a,b]$, is the Jordan variation of $f$ over $[a,b]$. 

By $AC_{loc}(\R)$ we denote the space of the functions that are locally absolutely continuous, namely absolutely continuous on every interval $[a,b]\subset \R$. Finally, we put $AC(\R):= BV(\R)\cap AC_{loc}(\R)$ (\cite{AV-2015}). We recall that a function of bounded variation is a.e. differentiable and 
$$
\int_{\R} |f'(t)|\,dt \le V[f].
$$
Moreover, if $f\in AC(\R)$, then 
\begin{equation}\label{int-rep}
\int_{\R} |f'(t)|\,dt = V[f],
\end{equation}
namely there holds an integral representation for the variation of $f$.

Let us consider the following family of discrete operators, known as {\it generalized sampling series} (\cite{RIST1,BUFIST}):
$$
(S_w f)(t):=\sum_{k\in\Z} f\left({k\over w}\right)\chi(wt-k), \ \ t\in\R,\ w>0,
$$
where $f:\R\rightarrow\R$ and $\chi:\R\rightarrow\R$ is a kernel that satisfies the following assumptions:
\begin{description} 
\item{$(\chi_1)$}
$\chi \in L^1(\R)$ is continuous on $\R$ 
such that $\sum_{k\in\Z} \chi(u-k)=1$, $\forall u\in\R$;
\item{$(\chi_2)$} 
$A_{\chi}:=\sup_{u\in\R}\sum_{k\in\Z} |\chi(u-k) |<+\infty$, where the convergence of the series is uniform on the compact sets of $\R$. 
\end{description}
The above assumptions are quite standard when we deal with discrete families of approximation operators: see e.g., \cite{BBSV-06,BCS1,COVI1,COVI2,ORTA1}.

We point out that, with such assumptions, $(S_w f)$ are well-defined if, for example, $f\in BV(\R)$. Indeed in this case $f$ is bounded and hence, if $M>0$ is such that $|f(x)|\le M$, for every $x\in\R$, 
$$
|(S_w f)(t)| \le M \sum_{k\in\Z} |\chi(wt-k)| \le A_{\chi}M <+\infty, \ \ t\in\R,\ w>0,
$$
by $(\chi_2)$. 

Such operators are well known and widely studied in approximation theory (see, e.g. \cite{BUNE,RIST1,BUFIST,K-2011}). 

In particular we will study the problem of the convergence in variation for the generalized sampling series in the case of an {\it averaged kernel}, namely of the form
\begin{equation}\label{av_ker}
\bar\chi_m(t):={1\over m}  \int_{-{m\over 2}}^{m\over 2} \chi(t+v)\,dv,
\end{equation}
for some $m\in\N$, where $\chi:\R\longrightarrow\R$ is a kernel. 

It can be proved that, if $\chi$ satisfies ($\chi_1$) and ($\chi_2$), the corresponding averaged kernel $\bar\chi_m$ turns out to be differentiable and satisfies the same conditions. 
Indeed, 
$$
\| \bar\chi_m \|_1 \leq m^{-1} \int^{m/2}_{-m/2}\left( \int_{\R}|\chi(v+t)|\, dt\right)dv = \| \chi\|_1<+\infty,
$$ 
and 
$$
\sum_{k \in \Z}\bar\chi_m(u-k) = m^{-1} \int_{-m/2}^{m/2} \left( \sum_{k \in \Z} \chi(u+v-k)\right)\, dv\ =1.
$$
Similarly, $(\chi_2)$ can be proved.

The sampling series corresponding to $\bar\chi_m$ are therefore of the form
$$
(\bar S^m_w f)(t):=\sum_{k\in\Z} f\left({k\over w}\right)\bar\chi_m(wt-k), \ \ t\in\R,\ w>0.
$$
Let us point out that there are several examples of well-known kernels in the literature of approximation theory that are of the form (\ref{av_ker}): in Section \ref{sec-ex} we will present some of them.  

Since we can obviously write 
\begin{equation}\label{av_ker_der}
\bar\chi'_m(t):={1\over m}  \left[\chi\left(t+{m\over 2}\right)-\chi\left(t-{m\over 2}\right)\right], \ \ t\in\R,
\end{equation}
for $f\in BV(\R)$ the derivative of the generalized sampling series $\bar S^m_w f$ can be written in the form 
\begin{equation}\label{sampl_der}
(\bar S^m_w f)'(t)={w\over m} \sum_{k\in\Z} f\left({k\over w}\right)\left[\chi\left(wt-k+{m\over 2}\right)-\chi\left(wt-k-{m\over 2}\right)\right]:
\end{equation}
notice that, for $f\in BV(\R)$, such derivative exists for every $t\in\R$ since 
\begin{equation*}
|(\bar S^m_w f)'(t)|\le {w\over m} \left(\left|S_w f\left(t+{m\over 2w}\right)\right|+\left|S_w f\left(t-{m\over 2w}\right)\right|\right),
\end{equation*}
for every $w>0$ and $m\in\N.$ 

In particular, in the case $m=1$, we will write 
\begin{equation*}\label{av_ker_1}
\bar\chi(t):= \int_{-{1\over 2}}^{1\over 2} \chi(t+v)\,dv,
\end{equation*}
and
\begin{equation*}\label{sampl_der_1}
(\bar S_w f)'(t)=w \sum_{k\in\Z} f\left({k\over w}\right)\left(\chi\left(wt-k+{1\over 2}\right)-\chi\left(wt-k-{1\over 2}\right)\right), 
\end{equation*}
for every $t\in\R$ and $w>0$.

\vskip0.3cm

One of the main goals of the present paper will be to establish a relation between the derivative of the generalized sampling series $(\bar S^m_w f)$ and the sampling-Kantorovich operators associated to the derivative of $f$. We recall that the sampling-Kantorovich operators (\cite{BBSV-06,CLCOMIVI2,COMIVI1,COGA2}) are defined as
\begin{equation}\label{sKant}
(K_w f)(t):=\sum_{k\in\Z} w\left(\int_{k\over w}^{k+1\over w} f(u)\,du\right) \chi(wt-k),\ t\in\R.
\end{equation}
Notice that, assuming ($\chi_2$), the operators $K_w f$ are well-defined for $f\in BV(\R)$ since, as before, 
$$
|(K_w f)(t)|\le M \sum_{k\in\Z} |\chi(wt-k)| \le MA_{\chi} <+\infty,\ t\in\R, \ w>0.
$$
Moreover, if $f\in AC(\R)$,
\begin{align*}
|(K_w f')(t)| &\le \sum_{k\in\Z} w\left|\int_{k\over w}^{k+1\over w} f'(u)\,du\right| | \chi(wt-k)| \\ &= \sum_{k\in\Z} w\left|f\left({k+1\over w}\right)-f\left({k\over w}\right)\right| | \chi(wt-k)| 
\\&
\le w V[f] \sum_{k\in\Z} | \chi(wt-k)| \le w A_{\chi} V[f] <+\infty,\ \ t\in\R, \ w>0.
\end{align*}


\section{Main results}\label{sec-mres}

\begin{prop} \label{var-dim}
Let $f\in BV(\R)$. Then $\bar S^m_w f\in AC(\R)$, for every $w>0$, $m\in\N$ and 
$$
V[\bar S^m_w f]\le {1\over m} \Vert \chi\Vert_1 V[f].
$$
\end{prop}

\noindent {\bf Proof.} Let $w>0$ and $m \in \N$. By (\ref{sampl_der}), since $f\in BV(\R)$, the derivative of $\bar S^m_w f$ exists everywhere in $\R$ and 
\begin{equation*}
|(\bar S^m_w f)'(t)|\le {w\over m} \sum_{k\in\Z} \left|f\left({k\over w}\right)\right|\left(\left|\chi\left(wt-k+{m\over 2}\right)\right|+\left|\chi\left(wt-k-{m\over 2}\right)\right|\right).
\end{equation*}
Now, since $f\in BV(\R)$, there exists $M\in\R$ such that $|f(x)|\le M$, for every $x\in\R$ and so, by assumption ($\chi_2$),
\begin{equation*}
|(\bar S^m_w f)'(t)|\le {wM\over m} \sum_{k\in\Z} \left(\left|\chi\left(wt-k+{m\over 2}\right)\right|+\left|\chi\left(wt-k-{m\over 2}\right)\right|\right) \le {2wM\over m} A_{\chi}.
\end{equation*}
Thus, $(\bar S^m_w f)'$ is bounded, which implies that $ \bar S^m_w f\in AC_{loc}(\R)$.
We will now prove that $\bar S^m_w f\in BV(\R)$.

Since $ \bar S^m_w f\in AC_{loc}(\R)$, there holds 
\begin{equation}\label{rint}
V[\bar S^m_w f] =\sup_{[a,b]\subset\R}V_{[a,b]}[\bar S^m_w f]=\sup_{[a,b]\subset\R} \int_a^b |(\bar S^m_w f)'(t)|\,dt = \int_{\R} |(\bar S^m_w f)'(t)|\,dt.
\end{equation}
Moreover we can write 
\begin{align*}
(\bar S^m_w f)'(t) &={w\over m} \sum_{k\in\Z} f\left({k\over w}\right)\left(\chi\left(wt-k+{m\over 2}\right)-\chi\left(wt-k-{m\over 2}\right)\right) \\
&= {w\over m} \sum_{k\in\Z} f\left({k\over w}\right)\chi\left(wt-k+{m\over 2}\right)-{w\over m} \sum_{k\in\Z} f\left({k\over w}\right)\chi\left(wt-k-{m\over 2}\right)\\ &=:{w\over m}(S_1+S_2).
\end{align*}
Now, putting $\tilde k=k+m$ in the series $S_2$, 
$$
S_2= \sum_{\tilde k\in\Z} f\left({\tilde k-m\over w}\right)\chi\left(wt-\tilde k+{m\over 2}\right),
$$
and so
\begin{align*}
(\bar S^m_w f)'(t) 
&= {w\over m} \sum_{k\in\Z} \left[f\left({k\over w}\right)-f\left({k-m\over w}\right)\right] \chi\left(wt-k+{m\over 2}\right). 
\end{align*}
Therefore, by (\ref{rint}), 
\begin{align*}
V[\bar S^m_w f]&= \int_{\R}|(\bar S^m_w f)'(t) |\,dt \\
&\le {w\over m} \int_{\R} \sum_{k\in\Z} \left|f\left({k\over w}\right)-f\left({k-m\over w}\right)\right| \left|\chi\left(wt-k+{m\over 2}\right)\right|\,dt \\ &\le  {1\over m} \int_{\R}\sum_{k\in\Z} \left|f\left({k\over w}\right)-f\left({k-m\over w}\right)\right|\left|\chi\left(u\right)\right|\,du \le {1\over m}V[f] \Vert \chi\Vert_1.
\end{align*}
Therefore the proof is complete, taking into account that $f\in BV(\R)$ and $\chi \in L^1(\R)$.\hfill\qed

\vskip0.3cm

The inequality proved in the above proposition is a {\em variation detracting (or diminishing)-type property} (see e.g., \cite{BBST-03,Ag-2006,K-2011,KM-2014,AA-2016}). In particular, in case of non-negative kernels we have $\|\chi\|_1=1$ and we obtain the usual variation diminishing property. Indeed, since $\chi$ is continuous, as a consequence of the Poisson's summation formula we know that the assumption:
$$
\sum_{k \in \Z} \chi(u-k)\ =\ 1,
$$
for every $u \in \R$, is equivalent to:
\begin{equation}
\widehat{\chi}(2\pi k)\ =\ \left\{
\begin{array}{l}
1, \hskip1cm k=0,\\
0, \hskip1cm k\neq0,
\end{array}
\right.
\end{equation}
for $k \in \Z$, where $\widehat{\chi}(v)= \int_{\R}\chi(u)\, e^{-i u v}\, du$, denotes the Fourier transform of $\chi$. Thus, since $\chi$ is non-negative, we obtain:
$$
1\ =\ \widehat{\chi}(0)\ =\ \int_{\R}\chi(u)\, du\ =\ \|\chi\|_1.
$$
Moreover, in the general case of kernels with variable sign it is sufficient to take $m \in \N$ large enough, to obtain again the usual variation diminishing property.

The next Proposition establishes a relation between the derivative of the generalized sampling series $(\bar S^m_w f)$ and the sampling-Kantorovich operators associated to the derivative of $f$.

\begin{prop}\label{prop-der}
Let $f\in AC(\R),$ then for every $t\in \R$, 
$$
(\bar S^m_w f)'(t)\ =\ {1\over m} \sum_{i=1}^m \left(K_w f'\right)\left(t-{m-2(i-1)\over 2w}\right),
$$
$w>0$, $m\in\N$. 
\end{prop}

\noindent {\bf Proof.} By (\ref{sampl_der_1}), we have that    
\begin{align*}
(\bar S^m_w f)'(t) &={w \over m} \sum_{k\in\Z} f\left({k\over w}\right)\left[\chi\left(wt-k+{m\over 2}\right)-\chi\left(wt-k-{m\over 2}\right)\right] \\ 
&= {w \over m} \sum_{k\in\Z} \left[ \int_{0}^{k\over w} f'(u)\,du + f(0)\right]\left[\chi\left(wt-k+{m\over 2}\right)-\chi\left(wt-k-{m\over 2}\right)\right] \\ 
&=  {w \over m} \sum_{k\in\Z} \left[ \int_{0}^{k\over w} f'(u)\,du + f(0)\right]\chi\left(wt-k+{m\over 2}\right)  +\\ & -  {w\over m} \sum_{k\in\Z} \left[ \int_{0}^{k\over w} f'(u)\,du + f(0)\right]\chi\left(wt-k-{m\over 2}\right).
\end{align*}
Now let us put, in the first series, $\tilde k=k-m$: then
\begin{align*}
(\bar S^m_w f)'(t) &=   {w \over m} \sum_{\tilde k\in\Z} \left[ \int_{0}^{\tilde k+m\over w} f'(u)\,du + f(0)\right]\chi\left(wt-\tilde k-{m\over 2}\right)  +\\ & -  {w \over m} \sum_{k\in\Z} \left[ \int_{0}^{k\over w} f'(u)\,du + f(0)\right]\chi\left(wt-k-{m\over 2}\right)\\  & =  {w \over m} \sum_{k\in\Z}  \left(\int_{k\over w}^{k+m\over w} f'(u)\,du \right)\chi\left(wt- k-{m\over 2}\right)  \\ &= {w \over m} \sum_{k\in\Z} \left[ \left( \int_{k\over w}^{k+1\over w}  + \ldots + \int_{k+(m-1)\over w}^{k+m\over w} \right) f'(u)\,du \right]\chi\left(wt- k-{m\over 2}\right)
\\ &={w \over m} \sum_{k\in\Z} \left( \int_{k\over w}^{k+1\over w} f'(u)\,du \right) \chi\left(wt- k-{m\over 2}\right) + \ldots + \\ &+ {w \over m} \sum_{k\in\Z} \left(\int_{k+(m-1)\over w}^{k+m\over w} f'(u)\,du \right) \chi\left(wt- (k+(m-1))-{m - 2(m-1)\over 2}\right)
\\ &= {1\over m}  (K_w f')\left(t-{m\over 2w}\right) + \ldots + {1\over m} (K_w f')\left(t-{m-2(m-1)\over 2w}\right)
\\ &= {1\over m} \sum_{i=1}^m \left(K_w f'\right)\left(t-{m-2(i-1)\over 2w}\right),
\end{align*}
for every $t\in\R$. \hfill\qed 

\vskip0.3cm

Note that the above property is analogous to a well-known relation occurring between the Bernstein polynomials and their Kantorovich-type version: see e.g., \cite{A-2001,BBST-03}.

We are now ready to prove the main result about a characterization of the absolute continuity in terms of the convergence in variation for the generalized sampling series $(\bar S^m_w f)$.

\begin{thm}\label{th-conv}
Let $f\in BV(\R)$. Then $\displaystyle\lim_{w\to +\infty} V[\bar S^m_w f-f]=0$, $m\in\N$, if and only if $f\in AC(\R)$.
\end{thm}

\noindent {\bf Proof.} We firstly consider $f\in AC(\R)$. By Proposition \ref{var-dim}, $\bar S^m_w f\in AC(\R)$, for every $w>0$, $m\in\N$, and hence also $(\bar S^m_w f -f) \in AC(\R)$. Therefore, by (\ref{int-rep}) and Proposition \ref{prop-der},
\begin{align*}
V[\bar S^m_w f-f] &=\int_{\R} |(\bar S^m_w f-f)'(t)|\,dt = \int_{\R} |(\bar S^m_w f)'-f'(t)|\,dt \\
&= \int_{\R} \left|{1\over m}\, \sum_{i=1}^m \left(K_w f'\right)\left(t-{m-2(i-1)\over 2w}\right) - {1\over m}\, m f'(t)\right|\,dt \\ 
& \le {1\over m} \int_{\R} \left| \sum_{i=1}^m \left(K_w f'\right)\left(t-{m-2(i-1)\over 2w}\right) - \sum_{i=1}^m f' \left(t-{m-2(i-1)\over 2w}\right) \right|\, dt\\
& +\ {1\over m} \int_{\R} \left| \sum_{i=1}^m f' \left(t-{m-2(i-1)\over 2w}\right) - m\, f'(t) \right|\, dt\\
& \le {1\over m} \sum_{i=1}^m \int_{\R}\left|\left(K_w f'\right)\left(t-{m-2(i-1)\over 2w}\right) - f' \left(t-{m-2(i-1)\over 2w}\right) \right|\, dt\\
& + {1\over m} \sum_{i=1}^m \int_{\R} \left| f' \left(t-{m-2(i-1)\over 2w}\right) - f'(t) \right|\, dt\\
& =\ \int_{\R}\left|\left(K_w f'\right)\left(t\right) - f' \left(t\right) \right|\, dt + {1\over m} \sum_{i=1}^m \int_{\R} \left| f' \left(t-{m-2(i-1)\over 2w}\right) - f'(t) \right|\, dt
 \\ &=:J\ +\ {1 \over m}\, \left\{I_1+\ldots+I_m\right\}.
\end{align*}
About $J$, since obviously
$$
J= \Vert K_w f'-f'\Vert_1
$$
and $f'\in L^1(\R)$, by Corollary 5.2 of \cite{BBSV-06}
\footnote{Notice that all the assumptions on kernels for such result are satisfied
(see Lemma 3.1 and Remark 3.2 of \cite{BBSV-06}).} 
we have that $J \rightarrow 0$, as $w\rightarrow +\infty$. Moreover, 
$$
I_1= \left\Vert f'\left(\cdot -{m\over 2w}\right)-f'(\cdot)\right\Vert_1\rightarrow 0,
$$
as $w\rightarrow +\infty$, by the continuity in $L^1$ of the translation operator, and analogously 
$$
I_j =\left\Vert f'\left(\cdot -{m-2(j-1)\over 2w}\right)-f'(\cdot)\right\Vert_1\rightarrow 0,
$$ 
for every $j=2,\ldots m$, as $w\rightarrow +\infty$. Therefore the first part of the theorem is proved. 

For the converse implication, notice that $\bar S^m_w f\in AC(\R)$, by Proposition \ref{var-dim}. Therefore, if $\displaystyle\lim_{w\to +\infty} V[\bar S^m_w f-f]=0$, recalling that $AC(\R)$ is a closed subspace of $BV(\R)$ with respect to the topology induced by the semi-norm defined by the total variation $V[\cdot]$ (see e.g. \cite{BBST-03}), we conclude that $f\in AC(\R)$.
\hfill\qed  


\section{Examples of kernel functions}\label{sec-ex}

In the literature, there are several examples of kernels (according to the definition given in Section \ref{sec-not}) which are of the averaged-type (\ref{av_ker}).

As first example, we present the case generated by the well-known central B-spline of order $n \in \N^+$, defined as follows:
\begin{equation} \label{splines}
 M_n(x)\ :=\ \frac{1}{(n-1)!} \sum^n_{i=0}(-1)^i \binom{n}{i} 
       \left(\frac{n}{2} + x - i \right)^{n-1}_+,     \hskip0.5cm   x \in \R,
\end{equation}
where $(x)_+ := \max\left\{x,0 \right\}$ denotes ``the positive part'' of $x \in \R$ (see e.g., \cite{MA1,BUNE,UN1,SP1}).

\begin{figure}
\centering
\includegraphics[scale=0.3]{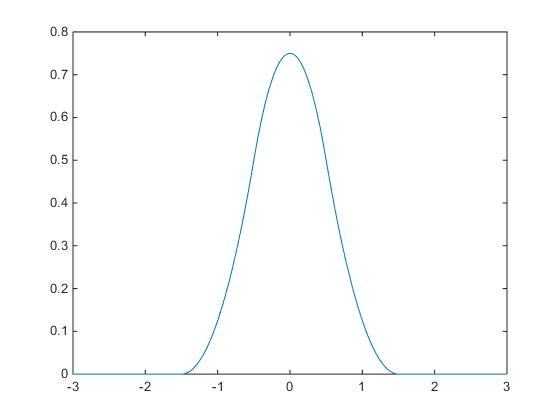}
\includegraphics[scale=0.3]{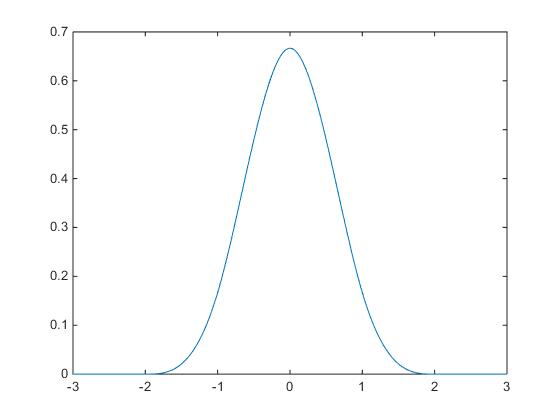}
\caption{The central B-spline of order $2$ (left), and its corresponding averaged kernel, i.e., the central B-spline of order $3$ (right).} \label{fig1}
\end{figure}

  The functions $M_n(x)$ (for some plots see Fig. \ref{fig1}) are non-negative, continuous with compact support contained in $[-n/2,n/2]$, and satisfy conditions $(\chi_1)$ and $(\chi_2)$.  
   
   In particular, the singularity assumption $\sum_{k \in \Z}M_n(u-k)=1$, for every $u \in \R$, follows as a consequence of the Poisson's summation formula, taking into account that $\widehat{M}_n$, i.e., the Fourier transform of $M_n$, is such that $\widehat{M}_n(2 \pi k)=0$, if $k \in \Z \setminus\left\{0\right\}$, and $\widehat{M}_n(0)=1$; for more details see e.g., \cite{BUNE}.

  Now, let us denote by
\begin{equation}
\bar M_{n,m}(t)\ :=\ m^{-1}\int_{-m/2}^{m/2}M_n(t+v)\, dv, \ \ t \in \R,
\end{equation}
the averaged B-spline kernel of order $n \in \N$. Recalling the following well-known property of the central B-spline, i.e., 
$$
M_n'(t)\ =\ M_{n-1}(t+1/2)\, -\, M_{n-1}(t-1/2), \ \ t \in \R, \ \ (n \geq 2)
$$
we have that, for $m=1$, 
$$
\bar M_{n,1}'(t)\ =\ M_{n}(t+1/2)\, -\, M_{n}(t-1/2)\ =\ M_{n+1}'(t), \ \ t \in \R, \ \  (n \geq 1),
$$
i.e., $\bar M_{n,1}(t)=M_{n+1}(t)+k$, $k \in \R$. Now, since $\bar M_{n,1}$ belongs to $L^1(\R)$, we must have $k=0$ and therefore we conclude that
\begin{equation}
\bar M_{n,1}(t)\ =\ M_{n+1}(t), \ \  t \in \R, 
\end{equation}
for every $n \in \N$, namely, the averaged kernel with $m=1$ generated by a central B-spline of order $n$ is a B-spline itself of order $n+1$ (see Fig. \ref{fig1} again). 
In view of what has been previously established, the following Corollary can be stated. 
\begin{cor} \label{cor1}
Let $f \in BV(\R)$. Denoting by $S^{M_n}_w$ the generalized sampling series based upon the central B-spline $M_n$ of order $n\ge 2$, we have that
$$
\lim_{w \to +\infty} V[S^{M_n}_w f - f]\ =\ 0
$$
if and only if $f \in AC(\R)$.
\end{cor}

\begin{figure}
\centering
\includegraphics[scale=0.3]{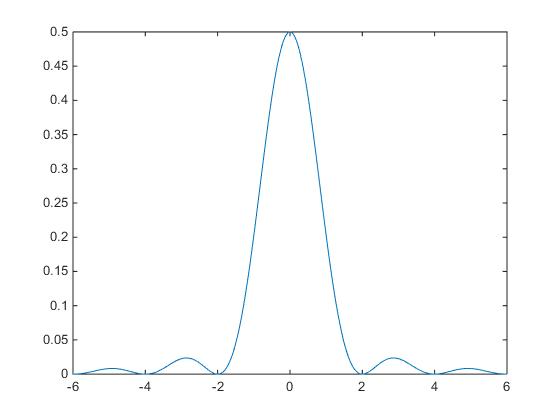}
\includegraphics[scale=0.3]{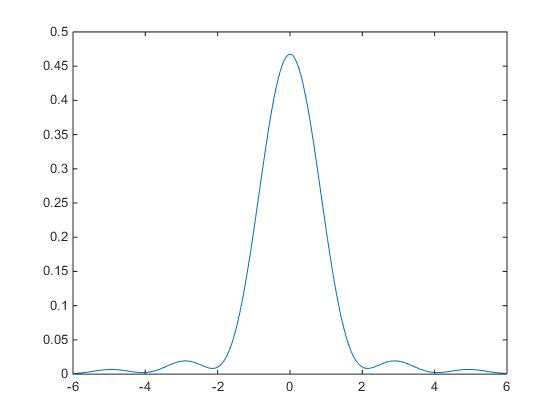}
\caption{The Fej\'er's kernel (left), and its corresponding averaged kernel with $m=1$ (right).} \label{fig2}
\end{figure}

In literature also examples of averaged kernels with unbounded support can be found. For instance, we can mention the Lanczos's kernel, defined by
$$
\bar \chi^s_m(t)\ :=\ m^{-1} \int_{-m/2}^{m/2}\mbox{sinc}(t+2v)\, dv, \ \  t \in \R,
$$
where
$$
\mbox{sinc}(t)\, :=\ \left\{
\begin{array}{l} 
\displaystyle {\sin{\pi t}\over\displaystyle \pi t}, \hskip0.5cm t \neq 0,\\
\\
1, \ \ t=0.
\end{array}
\right.
$$
Unfortunately, the $\mbox{sinc}$-function does not belong to $L^1(\R)$, then in this case the previous theory fails. However, this problem can be solved by considering a $\mbox{sinc}^2$-type kernel, such as the Fej\'er's kernel (see Fig. \ref{fig2}, left), defined by
$$
F(x)\ :=\ {1 \over 2}\,  \mbox{sinc}^2\left(\frac{x}{2}\right), \ \ x \in \R.
$$
The Fej\'er's kernel satisfies assumptions $(\chi_1)$ and $(\chi_2)$ (see e.g., \cite{COGA1,COVI4,COMIVI1}), and the corresponding averaged-kernel (see Fig. \ref{fig2}, right) takes now the form
$$
\bar F_m(t)\ :=\ m^{-1}\int_{-m/2}^{m/2}F(t+v)\, dv\ =\ {1 \over 2m}\int_{-m/2}^{m/2}\mbox{sinc}^2\left(\frac{t+v}{2}\right)\, dv, \ \ t\in \R.
$$
In particular, we can observe that $\bar F_m(t)$ turns out to be a bandlimited kernel, since $F(x)$ is bandlimited itself. The latter property can be viewed as a general fact; indeed, observing by (\ref{av_ker_der}) that $\bar \chi_m' \in L^1(\R)$, its Fourier transform can be computed and there holds
$$
\widehat{\bar \chi_m'}(v)\ =\ {\widehat{\chi}(v) \over m}\left[ e^{i m /2}\, -\,  e^{-i m /2} \right], \hskip1cm v \in \R.
$$
Now recalling that, in general, $\widehat{\bar \chi_m'}(v) = i v\, \widehat{\bar \chi_m}(v)$, $v \in \R$, it turns out that:
$$
\widehat{\bar \chi_m}(v)\ =\ m^{-1}\, {\widehat{\chi}(v)\over i v}\, \left[ e^{i m /2}\, -\,  e^{-i m /2} \right]\ =\ 2m^{-1}\, {\widehat{\chi}(v)\over v}\, \sin(m/2),
$$
$v \in \R\setminus\left\{ 0\right\}$. By the above equality, we conclude that $\bar \chi_m$ is bandlimited if and only if $\chi$ is bandlimited.

\begin{figure}
\centering
\includegraphics[scale=0.3]{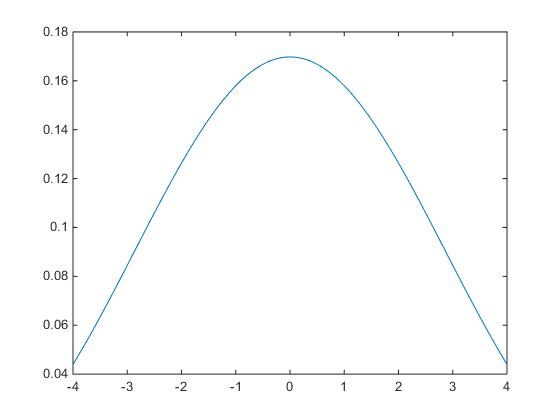}
\includegraphics[scale=0.3]{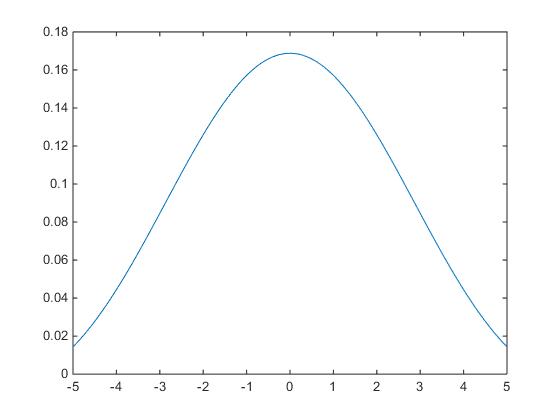}
\caption{The Bochner-Riesz kernel with $\gamma=2$ (left) and its corresponding averaged kernel with $m=1$ (right).} \label{fig3}
\end{figure}

  Other examples of bandlimited kernels are provided, e.g., by the Bochner-Riesz kernels (\cite{RS-2008}), defined by:
$$
b_{\gamma}(x)\ :=\ {2^{\gamma} \over \sqrt 2 \pi} \Gamma(\gamma+1)\, |x|^{-1/2 - \gamma}\, J_{1/2 + \gamma}(|x|), \hskip1cm x \in \R,
$$  
for $\gamma>0$, where $J_{\lambda}$ is the Bessel function of order $\lambda$ \cite{B-2012} (see Fig. \ref{fig3}, left). The corresponding averaged Bochner-Riesz kernels can be generated as follows:
$$
\bar b_{\gamma, m}(t)\ :=\ m^{-1} {2^{\gamma} \over \sqrt 2 \pi} \Gamma(\gamma+1) \int_{-m/2}^{m/2}|t+v|^{-1/2 - \gamma}\, J_{1/2 + \gamma}(|t+v|)\, dv, 
$$
$t \in \R$ (see Fig. \ref{fig3}, right). 
Now, denoting by $\bar S^m_w$ the generalized series based upon the kernel $\bar \chi_m$, with $\bar\chi_m(t)=\bar F_m(t)$ or $\bar\chi_m(t)=\bar b_{\gamma, m}(t)$, we can write what follows.
\begin{cor}
Let $f \in BV(\R)$, and $\bar\chi_m(t)=\bar F_m(t)$ or $\bar\chi_m(t)=\bar b_{\gamma, m}(t)$ be fixed, for some $m\in\N$. Then
$$
\lim_{w \to +\infty} V[\bar S^m_w f - f]\ =\ 0
$$
if and only if $f \in AC(\R)$.
\end{cor}

  By following the above procedure, several examples of kernels for which the previous theory holds can be given, see e.g., \cite{PZ-98,DOLU1,KT-2009,A-2013,AV-2013}.

\section{Final remarks and conclusions}

In this paper we prove a characterization of the absolute continuity in terms of the convergence in variation by means of the generalized sampling series $\bar S_w^m$. Such sampling series are based upon averaged kernels that do not need to be necessarily bandlimited. The crucial point of our approach is the possibility to establish a relation between the generalized sampling series and their Kantorovich-type version and the fact that the operators are based on averaged kernels. Actually, this is no restrictive, since there are many examples of kernels of averaged type well-known in approximation theory: among them, the central B-splines of order $n$.  As shown in Corollary \ref{cor1}, the generalized sampling series with averaged kernel generated by the central B-splines of order $n$ and $m=1$ coincide with the usual sampling series based upon the central B-splines of order $n+1$. It is well-known that the central B-splines are not bandlimited, therefore by the proposed approach we are able to treat a situation that was not covered by the convergence results proved in \cite{BV}, where the kernels are bandlimited and the function $f$ belongs to the Bernstein space (in general strictly contained in $AC(\R)$).  
Moreover, here we obtain not only a result of convergence of variation, but a complete characterization of $AC(\R)$ in terms of convergence in variation by means of the generalized sampling series, similarly to what happens, for example, working with the classical convolution integral operators.

\section*{Acknowledgments}

The authors are members of the Gruppo  
Nazionale per l'Analisi Matematica, la Probabilit\'a e le loro  
Applicazioni (GNAMPA) of the Istituto Nazionale di Alta Matematica (INdAM). 

\noindent The authors are partially supported by the "Department of Mathematics and Computer Science" of the University of Perugia (Italy). Moreover, D. Costarelli holds a research grant (Post-Doc) funded by the INdAM, and together with L. Angeloni, they have been partially supported within the 2017 GNAMPA-INdAM Project ``Approssimazione con operatori discreti e problemi di minimo per funzionali del calcolo delle variazioni con applicazioni all'imaging''.


\vskip0.1cm



\begin{thebibliography}{00}

\bibitem{AA-2016} U. Abel, O. Agratini, {\em On the variation detracting property of operators of Balazas and Szabados}, Acta Math. Hungar., {\bf 150} (2) (2016) 383--395.

\bibitem{A-2001} O. Agratini, {\em An approximation process of Kantorovich type}, Math. Notes, Miskolc, {\bf 2} (1), (2001), 3-10.

\bibitem{Ag-2006} O. Agratini, {\em On the variation detracting property of a class of operators}, Appl. Math. Lett., {\bf 19}(11) (2006), 1261--1264.

\bibitem{AD-2014} L. Ambrosio, S. Di Marino, {\em Equivalent definitions of BV space and of total variation on metric measure spaces}, Journal of Functional Analysis, {\bf 266} (7) (2014), 4150--4188.

\bibitem{A-2013} L. Angeloni, {\em Approximation results with respect to multidimensional $\varphi$-variation for nonlinear integral operators}, Z. Anal. Anwendungen, {\bf 32} (1) (2013), 103--128.

\bibitem{AV-2013} L. Angeloni, G. Vinti, {\em Approximation in variation by homothetic operators in multidimensional setting}, Differential Integral Equation, {\bf 26} (5-6) (2013), 655--674.

\bibitem{AV-2014} L. Angeloni, G. Vinti, {\em Convergence and rate of approximation in $BV^{\varphi}(\R_+^N)$ for a class of Mellin integral operators}, Atti della Accademia Nazionale dei Lincei, Classe di Scienze Fisiche, Matematiche e Naturali, Rendiconti Lincei Matematica e Applicazioni, {\bf 25} (3) (2014), 217--232

\bibitem{AV-2015} L. Angeloni, G. Vinti, {\em Convergence in variation and a characterization of  the absolute continuity}, Integral Transforms Spec. Funct., {\bf 26} (10) (2015), 829--844.

\bibitem{BBST-03} C. Bardaro, P.L. Butzer, R.L. Stens, G. Vinti, \emph{ Convergence in variation and rates of approximation for Bernstein-type polynomials and singular convolution integrals}, Analysis, {\bf 23} (2003) 299-340.

\bibitem{BBSV-06}
C. Bardaro, P.L. Butzer, R.L. Stens, G. Vinti, \emph{Kantorovich-Type Generalized Sampling Series in the Setting of Orlicz Spaces}, Sampling Theory in Signal and Image Processing, \textbf{6} (2007), 29--52.

\bibitem{BV} C. Bardaro, G. Vinti, {\em General convergence theorem with respect to Cesari variation and applications}, Nonlinear Anal., {\bf 22} (4) (1994), 505--518.

\bibitem{BCS1} A. Boccuto, D. Candeloro, and A. R. Sambucini,
\emph{Vitali-type theorems for filter convergence related to vector lattice-valued modulars and applications to stochastic
processes}, J. Math. Anal. Appl.  \textbf{419}(2) (2014), 818--838.

\bibitem{B-2012} F. Bowman, {\em Introduction to Bessel functions}, Courier Corporation (2012).

\bibitem{BUFIST} P.L. Butzer, A. Fisher, R.L. Stens, {\em Generalized sampling approximation of multivariate signals: theory and applications}, Note di Matematica, {\bf 1} (10) (1990), 173-191.

\bibitem{BUNE} P.L. Butzer, R.J. Nessel, {\em Fourier Analysis and Approximation I}, Academic Press, New York-London, 1971.

\bibitem{CLCOMIVI2} F. Cluni, D. Costarelli, A.M. Minotti, G. Vinti, {\em Applications of sampling Kantorovich operators to thermographic images for seismic engineering}, Journal of Computational Analysis and Applications, {\bf 19} (4) (2015) 602-617.

\bibitem{COGA2} L. Coroianu, S.G. Gal, $L^p$- approximation by truncated max-product sampling operators of Kantorovich-type based on Fejer kernel, in print in: Journal of Integral Equations and Applications, (2017).

\bibitem{COMIVI1} D. Costarelli, A.M. Minotti, G. Vinti, {\em Approximation of discontinuous signals by sampling Kantorovich series}, Journal of Mathematical Analysis and Applications, {\bf 450} (2) (2017) 1083-1103.

\bibitem{COVI1} D. Costarelli, G. Vinti, {\em Rate of approximation for multivariate sampling Kantorovich operators on some functions spaces}, J. Int. Eq. Appl., {\bf 26} (4) (2014), 455--481.

\bibitem{COVI2} D. Costarelli, G. Vinti, {\em Degree of approximation for nonlinear multivariate sampling Kantorovich operators on some functions spaces}, Numerical Functional Analysis and Optimization, {\bf 36} (8) (2015), 964--990. 

\bibitem{COVI4} D. Costarelli, G. Vinti, {\em Approximation by max-product neural network operators of Kantorovich type}, Results in Mathematics, {\bf 69} (3) (2016), 505--519. 

\bibitem{DOLU1} M.N. Do, Y.M. Lu, {\em A theory for sampling signals from a union of subspaces}, IEEE transactions on signal processing, {\bf 56} (6) (2008), 2334-2345. 

\bibitem{K-2011} A. Kivinukk, {\em On some Shannon sampling series with the variation detracting property}, Proc. of the 9th Intern. Conf. on Sampling Theory and Applications, Singapore (2011), 1--4.

\bibitem{KM-2014} A. Kivinukk, T. Metsmägi, {\em The variation detracting property of some Shannon sampling series and their derivatives}, Sampl. Theory Signal Image Process, {\bf 13} (2014), 189--206.

\bibitem{KT-2009} A. Kivinukk, G. Tamberg, {\em Interpolating generalized Shannon sampling operators, their norms and approximation properties}, Sampling Th. Signal Image Processing, {\bf 8} (2009), 77-95.

\bibitem{KR-2011} P. Krejci, T. Roche, {\em Lipschitz continuous data dependence of sweeping processes in BV spaces}, Discrete Contin. Dyn. Syst. Ser. B, {\bf 15} (2011) 637--650.

\bibitem{LS-2015} M. Laczkovich, V.T. S\'os, {\em Functions of Bounded Variation}, In: Real Analysis, Springer, New York, (2015) 399--406.

\bibitem{MA1} M.J. Marsden, {\em An identity for spline functions with applications to variation-diminishing spline approximation}, Journal of Approximation Theory, {\bf 3} (1) (1970), 7-49.

\bibitem{ORTA1} O. Orlova, G. Tamberg, {\em On approximation properties of generalized Kantorovich-type sampling operators}, Journal of Approximation Theory, {\bf 201} (2016), 73--86.

\bibitem{PZ-98} A. Piriou, X. Zeng, {\em On the rate of convergence of two Bernstein–B\'ezier type operators for bounded variation functions}, Journal of Approximation Theory, {\bf 95} (3) (1998). 369--387.

\bibitem{RIST1} S. Ries, R.L. Stens, {\em Approximation by generalized sampling series}, In: Constructive Theory of Functions'84, Sofia, (1984), 746--756.

\bibitem{RS-2008} K. Runovski, H.J. Schmeisser, {\em On approximation methods generated by Bochner-Riesz kernels}, Journal of Fourier Analysis and Applications, {\bf 14} (1) (2008) 16--38.

\bibitem{SP1} H. Speleers, {\em Multivariate normalized Powell-Sabin B-splines and quasi-interpolants}, Computer Aided Geometric Design, {\bf 30} (1) (2013), 2--19.

\bibitem{UN1} M.A. Unser, Ten good reasons for using spline wavelets, In: Optical Science, Engineering and Instrumentation'97. International Society for Optics and Photonics, (1997), 422-431. 

\end{thebibliography}
\end{document}